\documentclass[12pt]{article}

\usepackage[cmtip,arrow]{xy}
\usepackage{amsmath,amssymb,enumerate,pb-diagram,pb-xy}
\def\today{\number\day .\number\month .\number\year}
\usepackage[applemac]{inputenc}

\def \1{{\bf 1}}

\def \Ad{\operatorname{Ad}}
\def \al{\alpha}

\def \bs{\backslash}
\def \C{{\mathbb C}}
\def \CA{{\cal A}}
\def \CF{{\cal F}}

\def \CJ{{\cal J}}

\def \df{\ \stackrel{\mbox{\rm\tiny def}}{=}\ }

\def \e{\emph}
\def \eps{\varepsilon}
\def \Ext{\operatorname{Ext}}
\def \g{{\mathfrak g}}
\def \ga{\gamma}
\def \Ga{\Gamma}
\def \H{{\mathbb H}}
\def \Hom{{\rm Hom}}

\def \k{{\mathfrak k}}

\def \Mod{\operatorname{Mod}}

\def \mg{{\rm mg}}

\def \N{{\mathbb N}}
\def \p{{\mathfrak p}}
\def \pa{{\rm par}}

\def \ph{\varphi}

\def \prf{{\bf Proof: }}
\def \Q{{\mathbb Q}}
\def \qed{\ifmmode\eqno \square 
		\else\noproof\vskip 12pt plus 3pt minus 9pt \fi}
\def \noproof{{\unskip\nobreak\hfill\penalty50\hskip2em\hbox{}%
     \nobreak\hfill $\square$\parfillskip=0pt%
     \finalhyphendemerits=0\par}}

\def \sm{\smallsetminus}

\def \umg{{\rm umg}}
\def \vol{{\rm vol}}

\def \z{{\mathfrak z}}

\def \({\left(}
\def \){\right)}
\def \={{\ =\ }}

\newcommand{\tto}[1]{\stackrel{#1}{\longrightarrow}}

\newcommand{\norm}[1]{|\hspace{-1pt}| #1|\hspace{-1pt}|}

\newcommand{\ol}[1]{\overline{#1}}
\newcommand{\stack}[2]{\genfrac{}{}{0pt}{1}{#1}{#2}}

\newtheorem{theorem}{Theorem}[section]
\newtheorem{conjecture}[theorem]{Conjecture}
\newtheorem{lemma}[theorem]{Lemma}

\newtheorem{proposition}[theorem]{Proposition}

\begin{document}

\pagestyle{myheadings} \markright{HIGHER ORDER COHOMOLOGY OF ARITHMETIC GROUPS}

\title{Higher order cohomology of arithmetic groups}
\author{Anton Deitmar}
\date{}
\maketitle

{\bf Abstract:}
Higher order cohomology of arithmetic groups is expressed in terms of $(\g,K)$-cohomology.
Generalizing results of Borel, it is shown that the latter can be computed using functions of (uniform) moderate growth.
A higher order versions of Borel's conjecture is stated, asserting that the cohomology can be computed using automorphic forms.

$$ $$

\tableofcontents

\newpage
\section*{Introduction}

In \cite{ES} we have defined higher order group cohomology in the following general context:
Let $\Ga$ be a group and $\Sigma$ a normal subgroup.
For a ring $R$ we define a sequence of functors $H_q^0$ from the category of $R[\Ga]$-modules to the category of $R$-modules.
First, for an $R[\Ga]$-module $V$, one defines $H_1^0(\Ga,\Sigma,V)=H^0(\Ga,V)=V^\Ga$ as the fixed point module.
Inductively, $H_{q+1}^0(\Ga,\Sigma,V)$ is the module of all $v\in V$ such that $\sigma v=v$ for every $\sigma\in\Sigma$ and $\ga v-v$ is in $H_q^1(\Ga,\Sigma,V)$ for every $\ga\in\Ga$.
For every $q\ge 1$ the functor $H_q^0(\Ga,\Sigma,\cdot)$ is left-exact and we define the higher order group cohomology as the right derived functor
$$
H_q^p\= R^pH_q^0.
$$
In the case of a Fuchsian group the choice $\Sigma=\Ga_\pa=$ the subgroup generated by all parabolic elements, turned out to be the adequate choice for an Eichler-Shimura isomorphy result to hold, see \cite{ES}.
For general arithmetic groups $\Ga\subset G$, where $G$ is a reductive linear group over $\Q$, a replacement for the Eichler-Shimura isomorphism is the isomorphism to $(\g,K)$-cohomology,
$$
H^p(\Ga,E)\ \cong\ H_{\g,K}^p(C^\infty(\Ga\bs G)\otimes E),
$$
where $E$ is a finite dimensional representation of $G$.
In this paper we present a higher order analogue of this result, i.e.,
 we will show isomorphy of higher order cohomology to $(\g,K)$-cohomology,
$$
H_q^p(\Ga,,\Sigma,E)\ \cong\ H_{\g,K}^p(H_q^0(\Ga,\Sigma, C^\infty(G))\otimes E).
$$
We will prove higher order versions of results of Borel by which one can compute the cohomology using spaces of functions with growth restrictions.
We also state a higher order version of the Borel conjecture, proved by Franke \cite{Franke}, that the cohomology can be computed using automorphic forms.

Note that if $\Hom(\Ga,\C)=0$, then $H_q^p=H_1^p=H^p$ for every $q\ge 1$.
Consequently, in the case of arithmetic groups, higher order cohomology is of interested only for rank-one groups.

\section{General groups}
Let $R$ be a commutative ring with unit. Let $\Ga$ be a group and $\Sigma\subset \Ga$ a normal subgroup.
Let $I$ denote the augmentation ideal in the group algebra $A=R[\Ga]$.
Let $I_\Sigma$ denote the augmentation ideal of $R[\Sigma]$. As $\Sigma$ is normal in $\Ga$, the set $AI_\Sigma$ is a 2-sided ideal in $A$. 
For $q\ge 1$ consider the ideal
$$
J_q\df I^q+ R[\Ga]I_\Sigma.
$$
So in particular, for $\Sigma=\{ 1\}$ one has $J_q=I^q$.
On the other end, for $\Sigma=\Ga$ one gets $J_q=I$ for every $q\ge 1$.
For an $A$-module $V$ define
$$
H_q^p(\Ga,\Sigma,V)\=\Ext_A^p(A/J_q,V).
$$
This is the higher order cohomology of the module $V$, see \cite{ES}.
Note that in the case $q=1$, we get back the ordinary group cohomology, so
$$
H_1^p(\Ga,\Sigma,V)\= H^p(\Ga,V).
$$
For convenience, we will sometimes suppress the $\Sigma$ in the notation, so we simply write $H_q^p(\Ga,V)$ or even $H_q^p(V)$ for $H_q^p(\Ga,\Sigma,V)$.

For an $R$-module $M$ and a set $S$ we write $M^S$ for the $R$-module of all maps from $S$ to $M$.
Then $M^\emptyset$ is the trivial module $0$. 
Up to isomorphy, the module $M^S$ depends only on the cardinality of $S$.
It therefore makes sense to define $M^c$ for any cardinal number $c$ in this way.
Note that $J_q/J_{q+1}$ is a free $R$-module.
Define
$$
N_{\Ga,\Sigma}(q)\df \dim_R J_q/J_{q+1}.
$$
Then $N_{\Ga,\Sigma}(q)$ is a possibly infinite cardinal number.

\begin{lemma}\label{1.1}
\begin{enumerate}[\rm (a)]
\item For every $q\ge 1$ there is a natural exact sequence
\begin{multline*}
0\to H_q^0(\Ga,V)\to H_{q+1}^0(\Ga,V)\to H^0(\Ga,V)^{N_{\Ga,\Sigma}(q)}\to\\
\to H_q^1(\Ga,V)\to H_{q+1}^1(\Ga,V)\to H^1(\Ga,V)^{N_{\Ga,\Sigma}(q)}\to\dots\\
\dots\to H_q^p(\Ga,V)\to H_{q+1}^p(\Ga,V)\to H^p(\Ga,V)^{N_{\Ga,\Sigma}(q)}\to\dots
\end{multline*}
\item Suppose that for a given $p\ge 0$ one has $H^p(\Ga,V)=0$.
Then it follows $H_q^p(\Ga,V)=0$ for every $q\ge 1$.
In particular, if $V$ is acyclic as $\Ga$-module, then $H_q^p(\Ga,V)=0$ for all $p,q\ge 1$.
\end{enumerate}\end{lemma}

\prf
Consider the exact sequence
$$
0\to J_q/J_{q+1}\to A/J_{q+1}\to A/J_q\to 0.
$$
As an $A$-module, $J_q/J_{q+1}$ is isomorphic to a direct sum $\bigoplus_\al R_\al$ of copies of $R=A/I$.
So we conclude that for every $p\ge 0$,
$$
\Ext_A^p(J_q/J_{q+1},V)\ \cong\ \prod_\al\Ext_A^p(R,V)\ \cong\ H^p(\Ga,V)^{N_{\Ga,\Sigma}(q)}.
$$
The long exact $\Ext$-sequence induced by the above short sequence is
\begin{multline*}
0\to\Hom_A(A/J_q,V)\to\Hom_A(A/J_{q+1},V)\to \Hom_A(J_q/J_{q+1},V)\to\\
\to\Ext_A^1(A/J_q,V)\to\Ext_A^1(A/J_{q+1},V)\to \Ext_A^1(J_q/J_{q+1},V)\to\\
\to\Ext_A^2(A/J_q,V)\to\Ext_A^2(A/J_{q+1},V)\to \Ext_A^2(J_q/J_{q+1},V)\to\dots
\end{multline*}
This is the claim (a).
For (b) we proceed by induction on $q$.
For $q=1$ the claim follows from $H_1^p(\Ga,V)=H^p(\Ga,V)$.
Assume the claim proven for $q$ and $H^p(\Ga,V)=0$.
As part of the above exact sequence, we have the exactness of
$$
H_q^p(\Ga,V)\to H_{q+1}^p(\Ga,V)\to H^p(\Ga,V)^{N_{\Ga,\Sigma}(q)}.
$$
By assumption, we have $H^p(\Ga,V)^{N_{\Ga,\Sigma}(q)}=0$ and by induction hypothesis the module $H_q^p(\Ga,V)$ vanishes.
This implies $H_{q+1}^p(\Ga,V)=0$ as well.
\qed

\begin{lemma}[Cocycle representation]
The module ${\rm H}_q^1(\Ga,V)$ is naturally isomorphic to
$$
\Hom_A(J_q,V)/\al(V),
$$
where $\al:V\to \Hom_A(J_q,V)$ is given by $\al(v)(m)=mv$.
\end{lemma}

\prf
This is Lemma 1.3 of \cite{ES}.
\qed

\section{Higher order cohomology of sheaves}
Let $Y$ be a topological space which is path-connected and locally simply connected.
Let $C\to Y$ be a normal covering of $Y$.
Let $\Ga$ be the fundamental group of $Y$ and let $X\tto\pi Y$ be the universal covering.
The fundamental group $\Sigma$ of $C$ is a normal subgroup of $\Ga$.

For a sheaf $\CF$ on $Y$ define
$$
H_q^0(Y,C,\CF)\df H_q^0(\Ga,\Sigma,H^0(X,\pi^*\CF)).
$$
Let $\Mod(R)$ be the category of $R$-modules, let $\Mod_R(Y)$ be the category of sheaves of $R$-modules on $Y$, and let $\Mod_R(X)_\Ga$ be the category of sheaves over $X$ with an equivariant $\Ga$-action.
Then $H_q^0(Y,C,\cdot)$ is a left exact functor from $\Mod_R(Y)$ to $\Mod(R)$.
We denote its right derived functors by $H_q^p(Y,C,\cdot)$ for $p\ge 0$.

\begin{lemma} Assume that the universal cover $X$ is contractible.
\begin{enumerate}[\rm (a)]
\item For each $p\ge 0$ one has a natural isomorphism $H_1^p(Y,C,\CF)\ \cong\ H^p(Y,\CF)$.
\item If a sheaf $\CF$ is $H^0(Y,\cdot)$-acyclic, then it is $H_q^0(Y,C,\cdot)$-acyclic.
\end{enumerate}
\end{lemma}

Note that part (b) allows one to use flabby or fine resolutions to compute higher order cohomology.

\prf
We decompose the functor $H^0(Y,C,\cdot)$ into the functors
$$
\Mod_R(Y)\tto{\pi^*}\Mod_R(X)_\Ga\tto{H^0(X,\cdot)}\Mod(R[\Ga])\tto{H^0(\Ga,\Sigma,\cdot)}\Mod(R).
$$
The functor $\pi^*$ is exact and maps injectives to injectives.
We claim that $H^0(X,\cdot)$ has the same properties.
For the exactness, consider the commutative diagram
$$
\begin{diagram}
\node{\Mod_R(X)_\Ga}\arrow{e,t}{H^0}\arrow{s,r}f
	\node{\Mod(R[\Ga])}\arrow{s,r}f\\
\node{\Mod_R(X)}\arrow{e,t}{H^0}
	\node{\Mod(R),}
\end{diagram}
$$
where the vertical arrows are the forgetful functors.
As $X$ is contractible, the functor $H^0$ below is exact.
The forgetful functors have the property, that a sequence upstairs is exact if and only if its image downstairs is exact.
This implies that the above $H^0$ is exact.
It remains to show that $H^0$ maps injective objects to injective objects.
Let $\CJ\in\Mod_R(X)_\Ga$ be injective and consider a diagram with exact row in $\Mod(R[\Ga])$, 
$$
\begin{diagram}
\node{0}\arrow{e}
	\node{M}\arrow{e}\arrow{s,r}{\ph}
	\node{N}\\
\node[2]{H^0(X,\CJ).}
\end{diagram}
$$
The morphism $\ph$ gives rise to a morphism $\phi:M\times X\to \CJ$, where $M\times X$ stands for the constant sheaf with stalk $M$.
Note that $H^0(X,\phi)=\ph$.
As $\CJ$ is injective, there exists a morphism $\psi:N\times X\to \CJ$ making the diagram
$$
\begin{diagram}
\node{0}\arrow{e}
	\node{M\times X}\arrow{e}\arrow{s,r}{\phi}
	\node{N\times X}\arrow{sw,r}{\psi}\\
\node[2]{\CJ}
\end{diagram}
$$
commutative.
This diagram induces a corresponding diagram on the global sections, which implies that $H^0(X,\CJ)$ is indeed injective.

For a sheaf $\CF$ on $Y$ it follows that
$$
H^p(Y,\CF)\=R^p(H^0(Y,\CF))\=R^pH^0(\Ga,\Sigma,\CF)\circ H^0_\Ga\circ\pi^*\= H_1^p(Y,C,\CF).
$$
Now let $\CF$ be acyclic.
Then we conclude $H_1^p(\CF)=0$ for every $p\ge 1$, so the $\Ga$-module $V=H^0(X,\pi^*\CF)$ is $\Ga$-acyclic.
The claim follows from Lemma \ref{1.1}.
\qed

\section{Arithmetic groups}
Let $G$ be a semisimple Lie group with compact center and let $X=G/K$ be its symmetric space.
Let $\Ga\subset G$ be an arithmetic subgroup which is torsion-free, and let $\Sigma\subset \Ga$ be a normal subgroup.
Let $Y=\Ga\bs X$, then $\Ga$ is the fundamental group of the manifold $Y$, and the universal covering $X$ of $Y$ is contractible.
This means that we can apply the results of the last section.

\begin{theorem}
Let $(\sigma,E)$ be a finite dimensional representation of $G$.
There is a natural isomorphism
$$
H_q^p(\Ga,\Sigma,E)\ \cong\ H_{\g,K}^p(H_q^0(\Ga,\Sigma, C^\infty(G))\otimes E),
$$
where the right hand side is the $(\g,K)$-cohomology.
\end{theorem}

\prf
Let $\CF_E$ be the locally constant sheaf on $Y$ corresponding to $E$.
Let $\Omega_Y^p$ be the sheaf of complex valued $p$-differential forms on $Y$.
Then $\Omega_Y^p\otimes\CF_E$ is the sheaf of $\CF_E$-valued differential forms.
These form a fine resolution of $\CF_E$:
$$
0\to\CF_E\to\C^\infty\otimes\CF_E\tto{d\otimes 1}\Omega_Y^1\otimes\CF_E\to\dots
$$
Since $\pi^*\Omega_Y^\bullet=\Omega_X^\bullet$, 
we conclude that $H^p_q(\Ga,\Sigma,E)$ is the cohomology of the complex $H_q^0(\Ga,\Sigma,H^0(X,\Omega_X^\bullet\otimes E))$.
Let $\g$ and $\k$ be the Lie algebras of $G$ and $K$ respectively, and let $\g=\k\oplus\p$ be the Cartan decomposition.
Then $H^0(X,\Omega^p\otimes\CF_E)\=(C^\infty(G)\otimes\bigwedge^p\p)^K\otimes E$. 
Mapping a form $\omega$ in this space to $(1\otimes x^{-1})\omega(x)$ one gets an isomorphism to $(C^\infty(G)\otimes\bigwedge^p\p\otimes E)^K$, where $K$ acts diagonally on all factors and $\Ga$ now acts on $C^\infty(G)$ alone.
The claim follows.
\qed

Let $U(\g)$ act on $C^\infty$ as algebra of left invariant differential operators.
Let $\norm\cdot$ be a norm on $G$, see \cite{Wall}, Section 2.A.2.
Recall that a function $f\in C^\infty(G)$ is said to be \e{of moderate growth}, if for every $D\in U(\g)$ one has $Df(x)=O(\norm x^a)$ for some $a>0$.
The function $f$ is said to be of \e{uniform moderate growth}, if the exponent $a$ above can be chosen independent of $D$,
Let $C_\mg^\infty(G)$ and $C_\umg^\infty(G)$ denote the spaces of functions of moderate growth and uniform moderate growth respectively.

Let $\z$ be the center of the algebra $U(\g)$.
Let $\CA(G)$ denote the space of functions $f\in C^\infty(G)$ such that
\begin{itemize}
\item $f$ is of moderate growth,
\item $f$ is right $K$-finite, and
\item $f$ is $\z$-finite.
\end{itemize}

\begin{proposition}\label{2.2}
\begin{enumerate}[\rm (a)]
\item For $\Omega=C_\umg^\infty(G),C_\mg^\infty(G),C^\infty(G)$ one has
$$
H_q^1(\Ga,\Sigma,\Omega)\= 0
$$
for every $q\ge 1$.
\item If $\Hom(\Ga,\C)\ne 0$, then one has
$$
H^1(\Ga,\CA(G))\ \ne\ 0.
$$
\end{enumerate}
\end{proposition}

\prf
In order to prove (a), it suffices by Lemma \ref{1.1} (b), to consider the case $q=1$.
A 1-cocycle is a map $\al:\Ga\to \Omega$ such that $\al(\ga\tau)=\ga\al(\tau)+\al(\ga)$ holds for all $\ga,\tau\in\Ga$.
We have to show that for any given such map $\al$ there exists $f\in \Omega$ such that $\al(\tau)=\tau f-f$.
To this end 
consider the symmetric space $X=G/K$ of $G$.
Let $d(xK,yK)$ for $x,y\in G$ denote the distance in $X$ induced by the $G$-invariant Riemannian metric.
For $x\in G$ we also write $d(x)=d(xK,eK)$.
Then the functions $\log\norm x$ and $d(x)$ are equivalent in the sense that there exists a constant $C>1$ such that
$$
\frac 1C d(x)\ \le\ \log\norm x\ \le\ Cd(x)
$$
or
$$
\norm x\ \le\ e^{Cd(x)}\ \le\ \norm x^{C^2}
$$
holds for every $x\in G$.
We 
define
$$
\CF\=\{ y\in G: d(y)< d(\ga y)\ \forall \ga\in\Ga\sm\{ e\}\}.
$$
As $\Ga$ is torsion-free, this is a fundamental domain for the $\Ga$ left translation action on $G$.
In other words, $\CF$ is open, its boundary is of measure zero, and there exists a set of representatives $R\subset G$ for the $\Ga$-action such that $\CF\subset R\subset \ol\CF$.
Next let $\ph\in C_c^\infty(G)$ with $\ph\ge 0$ and $\int_G\ph(x)\,dx=1$.
Then set $u=\1_\CF *\ph$, where $\1_A$ is the characteristic function of the set $A$ and $*$ is the convolution product $f*g(x)=\int_Gf(y)g(y^{-1}x)\,dy$.
Let $C$ be the support of $\ph$, then the support of $u$ is a subset of $\ol\CF C$ and the sum $\sum_{\tau\in\Ga} u(\tau^{-1}x)$ is locally finite in $x$.
More sharply, for a given compact unit-neighborhood $V$ there exists $N\in\N$ such that for every $x\in G$ one has
$$
\# \{ \tau\in\Ga: u(\tau^{-1} xV)\nsubseteq \{ 0\}\}\ \le \ N.
$$
This is to say, the sum is uniformly locally finite.
For a function $h$ on $G$ and $x,y\in G$ we write $L_yh(x)=h(y^{-1} x)$.
Then for a convolution product one has $L_y(f*g)=(L_yf)*g$, and so
$$
\sum_{\tau\in\Ga}u(\tau^{-1}x)\=\(\sum_{\tau\in\Ga}L_\tau\1_\CF\)*\ph.
$$
The sum in parenthesis is equal to one on the complement of a nullset.
Therefore,
$$
\sum_{\tau\in\Ga}u(\tau^{-1} x)\ \equiv\ 1.
$$
Set
$$
f(x)\= -\sum_{\tau\in\Ga}\al(\tau)(x)\,u(\tau^{-1} x).
$$
\begin{lemma}
The function $f$ lies in the space $\Omega$.
\end{lemma}

\prf
Since the sum is uniformly locally finite, it suffices to show that for each $\tau\in\Ga$ we have $\al(\tau)(x)u(\tau^{-1} x)\in \Omega$ where the $O(\norm\cdot^d)$ estimate is uniform in $\tau$.
By the Leibniz-rule it suffices to show this separately for the two factors $\al(\tau)$ and $L_\tau u$.
For $D\in U(\g)$ we have
$$
D(L_\tau u)=(L_\tau\1_\CF)*(D\ph).
$$
This function is bounded uniformly in $\tau$, hence $L_\tau u\in C_\umg^\infty(G)$.
Now $\al(\tau)\in\Omega$ by definition, but we need uniformity of growth in $\tau$.
We will treat the case $\Omega= C_\umg^\infty(G)$ here, the case $C_\mg^\infty$ is similar and the case $C^\infty(G)$ is trivial, as no growth bounds are required.

So let $\Omega= C_\umg^\infty(G)$ and set
$$
S\=\{ \ga\in\Ga\sm\{ e\} : \ga\ol\CF\cap\ol\CF\ne \emptyset\}.
$$
Then $S$ is a finite symmetric generating set for $\Ga$.
For $\ga\in\Ga$, let $\CF_\ga$ be the set of all $x\in G$ with $d(x) <d(\ga x)$.
Then
$$
\CF\=\bigcap_{\ga\in\Ga\sm\{ e\}}\CF_\ga
$$
Let $\tilde\CF\=\bigcap_{s\in S}\CF_s$. We claim that $\CF=\tilde \CF$.
As the intersection runs over fewer elements, one has $\CF\subset\tilde\CF$.
For the converse note that for every $s\in S$ the set $s\ol\CF/K$ lies in $X\sm\tilde\CF/K$, therefore $\CF/K$ is a connected component of $\tilde\CF/K$.
By the invariance of the metric, we conclude that $x\in\CF_\ga$ if and only if $d(xK,eK)<d(xK,\ga^{-1}K)$.
This implies that $\CF_\ga/K$ is a convex subset of $X$.
Any intersection of convex sets remains convex, therefore $\tilde\CF/K$ is convex and hence connected, and so $\tilde \CF/K=\CF/K$, which means $\tilde\CF=\CF$.

Likewise we get $\ol\CF\=\bigcap_{s\in S}\ol{\CF_s}$.
The latter implies that for each $x\in G\sm\ol\CF$ there exists $s\in S$ such that $d(s^{-1}x)<d(x)$.
Iterating this and using the fact that the set of all $d(\ga x)$ for $\ga\in\Ga$ is discrete, we find for each $x\in G\sm\ol\CF$ a chain of elements $s_1,\dots,s_n\in S$ 
such that $d(x) >d(s_1^{-1}x)>\dots>d(s_n^{-1}\cdots s_1^{-1} x)$
 and $s_n^{-1}\cdots s_1^{-1} x\in\ol\CF$.
The latter can be written as $x\in s_1\dots s_n\ol\CF$.
Now let $\tau\in\Ga$ and suppose $u(\tau^{-1} x)\ne 0$.
Then $x\in\ol\CF C$, so, choosing $C$ small enough, we can assume $x\in s\tau\ol\CF$ for some $s\in S\cap\{ e\}$.
As the other case is similar, we can assume $s=e$.
It suffices to assume $x\in\tau\CF$, as we only need the estimates on the dense open set $\Ga\CF$.
So then it follows $\tau= s_1\dots s_n$.

Let $D\in U(\g)$.
As $\al$ maps to $\Omega= C_\umg^\infty(G)$, for every $\ga\in\Ga$ there exist $C(D,\ga), a(\ga)>0$ such that
$$
|D\al(\ga)(x)|\ \le\ C(D,\ga)\,\norm x^{a(\ga)}.
$$
The cocycle relation of $\al$ implies
$$
\al(\tau)(x)\=\sum_{j=1}^n\al(\ga_j)(s_{j-1}^{-1}\dots s_1^{-1} x).
$$
We get
\begin{eqnarray*}
|D\al(\tau)(x)| &\le& \sum_{j=1}^n C(D,s_j)\norm{s_{j-1}^{-1}\dots s_1^{-1}x}^{a(s_j)}\\
&\le& \sum_{j=1}^n C(D,s_j)\, e^{Cd(s_{j-1}^{-1}\dots s_1^{-1}x)a(s_j)}\\
&\le& \sum_{j=1}^n C(D,s_j)\, e^{Cd(x)a(s_j)}\\
&\le& \sum_{j=1}^n C(D,s_j)\norm{x}^{C^2a(s_j)}\\
&\le& nC_0(D)\norm x^{a_0},
\end{eqnarray*}
where $C(D)=\max_j C(D,s_j)$ and $a_0\=C^2\max_jd(s_j)$.
It remains to show that $n$ only grows like a power of $\norm x$.
To this end let for $r>0$ denote $N(r)$ the number of $\ga\in\Ga$ with $d(\ga)\le r$.
Then a simple geometric argument shows that
$$
N(r)\=\frac 1{\vol\CF}\vol\(\bigcup_{\ga:d(\ga)\le r}\ga\CF/K\)\ \le\ C_1\vol (B_{2r}),
$$
where $B_{2r}$ is the ball of radius $2r$ around $eK$.
Note that for the homogeneous space $X$ there exists a constant $C_2>0$ such that $\vol B_{2r}\le e^{C_2 r}$.
Now $n\le N(d(x))$ and therefore
$$
n\ \le\ C_1\vol B_{2d(x)}\ \le\ C_1e^{C_2 d(x)}\ \le\ C_1 \norm x^{C_3}
$$
for some $C_3>0$.
Together it follows that there exists $C(D)>0$ and $a>0$ such that
$$
|D\al(\tau)(x)|\ \le\ C(D)\,\norm x^a.
$$
This is the desired estimate which shows that $f\in\Omega$.
The lemma is proven.
\qed

To finish the proof of part (a) of the proposition, we now compute for $\ga\in\Ga$,
\begin{eqnarray*}
\ga f(x)- f(x) &=& f(\ga^{-1} x)-f(x)\\
&=& \sum_{\tau\in\Ga}\al(\tau x) u(\tau^{-1} x) -
\al(\tau)(\ga^{-1} x)u(\tau^{-1}\ga^{-1} x)\\
&=& \sum_{\tau\in\Ga}\al(\tau) (x) u(\tau^{-1}x) 
+\al(\ga)(x)\sum_{\tau\in\Ga} u((\ga\tau)^{-1}x)\\
&& -\sum_{\tau\in\Ga}\al(\ga \tau)(x)u((\ga\tau)^{-1} x)
\end{eqnarray*}
The first and the last sum cancel and the middle sum is $\al(\ga)(x)$.
Therefore, part (a) of the  proposition is proven.

We now prove part (b).
Let $Q=C^\infty(G)/\CA(G)$.
We have an exact sequence of $\Ga$-modules
$$
0\to\CA(G)\to C^\infty(G)\to Q\to 0.
$$
This results in the exact seuence
$$
0\to \CA(G)^\Ga\to C^\infty(\Ga\bs G)\tto\phi Q^\Ga\to H^1(\Ga,\CA(G))\to 0.
$$
The last zero comes by part (a) of the proposition.
We have to show that the map $\phi$ is not surjective.
So let $\chi:\Ga\to\C$ be a non-zero group homomorphism and let $u\in C^\infty(G)$ as above with $\sum_{\tau\in\Ga} u(\tau^{-1} x)=1$, and $u$ is supported in $\ol\CF C$ for a small unit-neighborhood $C$.
Set
$$
h(x)\=-\sum_{\tau\in\Ga}\chi(\tau) u(\tau^{-1} x).
$$
Then for every $\ga\in\Ga$ the function
$$
h(\ga^{-1} x)-h(x)\=\chi(\ga)
$$
is constant and hence
lies in $\CA(G)^\Ga$. This means that the class $[h]$ of $h$ in $Q$ lies in the $\Ga$-invariants $Q^\Ga$.
As $\chi\ne 0$,the function $f$ is not in $C^\infty(\Ga\bs G)$, and therefore $\phi$ is indeed not surjective.
\qed

\begin{proposition}\label{2.4}
For every $q\ge 1$ there is an exact sequence of continuous $G$-homomorphisms,
$$
0\to H_q^0(\Ga,\Sigma,C_*^\infty(G))\tto\phi  H_{q+1}^0(\Ga,\Sigma,C_*^\infty(G))
\tto\psi C_*^\infty(\Ga\bs G)^{N_{\Ga,\Sigma}(q)}\to 0,
$$
where $\phi$ is the inclusion map and $*$ can be $\emptyset, \umg$, or $\mg$.
\end{proposition}

\prf
This follows from Lemma \ref{1.1} together with Propostion \ref{2.2} (a).
\qed

The space $C^\infty(G)$ carries a natural topology which makes it a nuclear topological vector space.
For every $q\ge 1$, the space $H_q^0(\Ga,\Sigma,C^\infty(G))$ is a closed subspace.
If $\Ga$ is cocompact, then one has the isotypical decomposition
$$
H_1^0(\Ga,\Sigma, C^\infty(G))\= C^\infty(\Ga\bs G)\= \ol{\bigoplus_{\pi\in\hat G}} C^\infty(\Ga\bs G)(\pi),
$$
and $C^\infty(\Ga\bs G)(\pi)\cong m_\Ga(\pi)\pi^\infty$,
where the sum runs over the unitary dual $\hat G$ of $G$, and for $\pi\in\hat G$ we write $\pi^\infty$ for the space of smooth vectors in $\pi$.
The multiplicity $m_\Ga(\pi)\in\N_0$ is the multiplicity of $\pi$ as a subrepresentation of $L^2(\Ga\bs G)$, i.e.,
$$
m_\Ga(\pi)\=\dim\Hom_G\(\pi,L^2(\Ga\bs G)\).
$$
Finally, the direct sum $\ol{\bigoplus}$ means the closure of the algebraic direct sum in $C^\infty(G)$.
We write $\hat G_\Ga$ for the set of all $\pi\in\hat G$ with $m_\Ga(\pi)\ne 0$.

Let $\pi\in\hat G$.
A smooth representation $(\beta,V_\beta)$ of $G$ is said to be \e{of type $\pi$}, if it is of finite length and every irreducible subquotient is isomorphic to $\pi^\infty$.
For a smooth representation $(\eta,V_\eta)$ we define the \e{$\pi$-isotype} as
$$
V_\eta(\pi)\df \ol{\bigoplus_{\stack{V_\beta\subset V_\eta}{\beta{\rm\ of\ type\ }\pi}}}V_\beta,
$$
where the sum runs over all subrepresentations $V_\beta$ of type $\pi$.

\begin{theorem}
Suppose $\Ga$ is cocompact and let $*\in\{ \emptyset,\mg,\umg\}$.
We write $V_q=V$.
For every $q\ge 1$ there is an isotypical decomposition
$$
V_q\=\ol{\bigoplus_{\pi\in\hat G_\Ga}}V_q(\pi),
$$
and each $V_q(\pi)$ is of type $\pi$ itself.
The exact sequence of Proposition \ref{2.4} induces an exact sequence
$$
0\to V_q(\pi)\to V_{q+1}(\pi)\to (\pi^\infty)^{m_\Ga(\pi)N_{\Ga,\Sigma}(q)}\to 0
$$
for every $\pi\in\hat G_\Ga$.
\end{theorem}

\prf
We will prove the theorem by reducing to a finite dimensional situation by means of considering infinitesimal characters and $K$-types.
For this let $\hat\z=\Hom(\z,\C)$ be the set of all algebra homomorphisms from $\z$ to $\C$.
For a $\z$-module $V$ and $\chi\in\hat\z$ let
$$
V(\chi)\df \{ v\in V: \forall z\in\z\ \exists n\in\N\ (z-\chi(z))^nv=0\}
$$
be the \e{generalized $\chi$-eigenspace}.
Since $\z$ is finitely generated, one has
$$
V(\chi)\= \{ v\in V: \exists n\in\N\ \forall z\in\z\ (z-\chi(z))^nv=0\}.
$$
For $\chi\ne \chi'$ in $\z$ one has $V(\chi)\cap V(\chi')=0$.
Recall that the algebra $\z$ is free in $r$ generators, where $r$ is the absolute rank of $G$.
Fix a set of generators $z_1,\dots, z_r$.
The map $\chi\mapsto (\chi(z_1),\dots,\chi(z_r))$ is a bijection $\hat\z\to\C^r$.
We equip $\hat\z$ with the topology of $\C^r$.
This topology does not depend on the choice of the generators $z_1,\dots,z_r$.

Let $\Ga\subset G$ be a discrete cocompact subgroup. 
Let $\hat\z_\Ga$ be the set of all $\chi\in\hat\z$ such that the generalized eigenspace $C^\infty(\Ga\bs G)(\chi)$ is non-zero.
The $\hat\z_\Ga$ is discrete in $\hat\z$, more sharply there exists 
$\eps_\Ga>0$ such that for any two $\chi\ne\chi'$ in $\hat\z_\Ga$ 
there is $j\in\{ 1,\dots,r\}$ such that $|\chi(z_j)-\chi'(z_j)|>\eps_\Ga$.

\begin{proposition}
Let $*\in\{ \emptyset,\mg,\umg\}$.
For every $q\ge 1$ and every $\chi\in\hat\z$ the space 
$V_q(\chi)=H_q^0(\Ga,\Sigma,C^\infty_*(G))(\chi)$ coincides with
$$
\bigcap_{z\in\z}\ker(z-\chi(z))^{2^{q-1}},
$$
and is therefore a closed subspace of $V_q$.
The representation of $G$ on $V_q(\chi)$ is of finite length.

The space $V_q(\chi)$ is non-zero only if $\chi\in\hat\z_\Ga$.
One has a decomposition
$$
H_q^0(\Ga,\Sigma,C^\infty_*(G))\=\ol{\bigoplus_{\chi\in\hat\z_\Ga}}H_q^0(\Ga,\Sigma,C^\infty_*(G))(\chi).
$$
The exact sequence of Proposition \ref{2.4} induces an exact sequence
$$
0\to V_q(\chi)\to  V_{q+1}(\chi)
\to \bigoplus_{\pi\in\hat G_\chi}m_\Ga(\pi)N_{\Ga,\Sigma}(q)\pi
\to 0.
$$
\end{proposition}

\prf
All assertions, except for the exactness of the  sequence, are clear for $q=1$.
We proceed by induction.
Fix $\chi\in\hat\z_\Ga$.
Since $V_q(\chi)=V_q\cap V_{q+1}(\chi)$ one gets an exact sequence
$$
0\to V_q(\chi)\to V_{q+1}(\chi)\tto{\psi_\chi}V_1(\chi)^{N_{\Ga,\Sigma}(q)}.
$$
Let $v\in V_1(\chi)^{N_{\Ga,\Sigma}(q)}$.
As $\psi$ is surjective, one finds $u\in V_{q+1}$ with $\psi(u)=v$.
We have to show that one can choose $u$ to lie in $V_{q+1}(\chi)$.
We have $(z-\chi(z))v=0$ for every $z\in\z$.
Therefore $(z-\chi(z))u\in V_q$.
Inductively we assume the decomposition to holds for $V_q$, so we can write
$$
(z_j-\chi(z_j))u\=\sum_{\chi'\in\hat\z_\Ga}u_{j,\chi'},
$$
for $1\le j\le r$ and $u_{j,\chi'}\in\ker(z-\chi(z))^{2^{q-1}}$ for every $z\in\z$.
For every $\chi'\in\hat\z_\Ga\sm\{\chi\}$ we fix some index $1\le j(\chi')\le r$ with $|\chi(z_{j(\chi')}-\chi'(z_{j(\chi')}|>\eps_\Ga$.
On the space
$$
\ol{\bigoplus_{\chi':j(\chi')=j}}V_q(\chi')
$$
the operator $z_j-\chi(z_j)$ is invertible and the inverse $(z_j-\chi(z_j))^{-1}$ is continuous.
We can replace $u$ with
$$
u-\sum_{\chi'\in\hat\z_{\Ga}\sm\{\chi\}}(z_{j(\chi')}-\chi(z_{j(\chi')}))^{-1}u_{j(\chi'),\chi'}.
$$
We end up with $u$ satisfying $\psi(u)=v$ and
$$
(z_1-\chi(z_1))\cdots (z_r-\chi(z_r))u\ \in\ V_q(\chi)\=\bigcap_{z\in\z}\ker(z-\chi(z))^{2^{q-1}}.
$$
So for every $z\in\z$ one has
$$
0\= (z_1-\chi(z_1))\cdots (z_r-\chi(z_r))(z-\chi(z))^{2^{q-1}}u,
$$
which implies
$$
(z-\chi(z))^{2^{q-1}}u\in\ker\((z_1-\chi(z_1))\cdots (z_r-\chi(z_r))\).
$$
As the set $\z_\Ga$ is countable, one can, depending on $\chi$, choose the generators $z_1,\dots,z_r$ in a way that $\chi(z_j)\ne \chi'(z_j)$ holds for every $j$ and every $\chi'\in\hat\z_\Ga\sm\{\chi\}$.
Therefore the operator $(z_1-\chi(z_1))\cdots (z_r-\chi(z_r))$ 
is invertible on $V_{q}(\chi')$ for every 
$\chi'\in\hat\z_\Ga\sm\{\chi\}$ and so it follows 
$(z-\chi(z))^{2^{q-1}}u\in V_q(\chi)\subset \ker(z-\chi(z))^{2^{q-1}}$ 
and therefore $u\in\ker((z-\chi(z))^{2^q})$. 
Since this holds for every $z$ it follows $u\in V_{q+1}(\chi)$ and so $\psi_\chi$ is indeed surjective. One has an exact sequence
$$
0\to V_q(\chi)\to V_{q+1}(\chi)\to V_1(\chi)^{N_{\Ga,\Sigma}(q)}\to 0.
$$
Taking the sum over all $\chi\in\z_\Ga$ we arrive at an exact sequence
$$
0\to V_q\to \ol{\bigoplus_{\chi\in\hat\z_\Ga}}V_{q+1}(\chi)\to V_1^{N_{\Ga,\Sigma}(q)}\to 0.
$$
Hence we get a commutative diagram with exact rows
$$\divide\dgARROWLENGTH by2
\begin{diagram}
\node{0}\arrow{e}
	\node{V_q}\arrow{e}\arrow{s,r}{=}
	\node{\ol{\bigoplus_{\chi\in\hat\z_\Ga}}V_{q+1}(\chi)}\arrow{e}\arrow{s,r,J}{i}
	\node{V_1^{N_{\Ga,\Sigma}(q)}}\arrow{e}\arrow{s,r}{=}
	\node{0}\\
\node{0}\arrow{e}
	\node{V_q}\arrow{e}
	\node{V_{q+1}}\arrow{e}
	\node{V_1^{N_{\Ga,\Sigma}(q)}}\arrow{e}
	\node{0,}	
\end{diagram}
$$
where $i$ is the inclusion.
By the 5-Lemma, $i$ must be a bijection.
The proposition follows.
\qed

We now finish the proof of the theorem.
We keep the notation $V_q$ for the space $\H_q^0(\Ga,\Sigma,C_*^\infty(G))$.
For a given $\chi\in \hat\z_\Ga$ the $G$-representation $V_q(\chi)$ is of finite length, so the $K$-isotypical decomposition
$$
V_q(\chi)\=\ol{\bigoplus_{\tau\in\hat K}}V_q(\chi)(\tau)
$$
has finite dimensional isotypes, i.e., $\dim V_q(\chi)(\tau)<\infty$.
Let $U(\g)^K$ be the algebra of all $D\in U(\g)$ such that $\Ad(k)D=D$ for every $k\in K$.
Then the action of $D\in U(\g)$ commutes with the action of each $k\in K$, and so $K\times U(\g)^K$ acts on every smooth $G$-module.
For $\pi\in\hat G$ the $K\times U(\g)^K$-module $V_\pi(\tau)$ is irreducible and $V_\pi(\tau)\cong V_{\pi'}(\tau')$ as a $K\times U(\g)^K$-module implies 
$\pi=\pi'$ and $\tau=\tau'$, see \cite{Wall}, Proposition 3.5.4.
As $V_q(\chi)(\tau)$ is finite dimensional. one gets
$$
V_q(\chi)(\tau)\=\bigoplus_{\stack{\pi\in\hat G}{\chi_\pi=\chi}}V_q(\chi)(\tau)(\pi),
$$
where $V_q(\chi)(\tau)(\pi)$ is the largest $K\times U(\g)^K$-submodule of $V_q(\chi)(\tau)$ with the property that every irreducible subquotient is isomorphic to $V_\pi(\tau)$.
Let 
$$
V_q(\pi)\=\ol{\bigoplus_{\tau\in\hat K}}V_q(\chi_\pi)(\tau)(\pi).
$$
The claims of the theorem follow from the proposition.
\qed

\section{The higher order Borel conjecture}

Let $(\sigma,E)$ be a finite dimensional representation of $G$. In \cite{Borel}, 
A. Borel has shown that the inclusions $C_\umg^\infty(G)\hookrightarrow C_\mg^\infty(G)\hookrightarrow C^\infty(G)$ induce isomorphisms in cohomology:
\begin{multline*}
H_{\g,K}^p(H^0(\Ga,C_\umg^\infty(G))\otimes E)\tto\cong 
H_{\g,K}^p(H^0(\Ga,C_\mg^\infty(G))\otimes E)\\
\tto\cong 
H_{\g,K}^p(H^0(\Ga,C^\infty(G))\otimes E).
\end{multline*}

In \cite{Franke}, J. Franke proved a conjecture of Borel stating that the inclusion $\CA(G)\hookrightarrow C^\infty(G)$ induces an isomorphism
$$
H_{\g,K}^p(H^0(\Ga,\CA(G))\otimes E)\tto\cong 
H_{\g,K}^p(H^0(\Ga,C^\infty(G))\otimes E).
$$

\begin{conjecture}
[Higher order Borel conjecture]
For every $q\ge 1$, the inclusion $\CA(G)\hookrightarrow C^\infty(G)$ induces an isomorphism
$$
H_{\g,K}^p(H_q^0(\Ga,\Sigma,\CA(G))\otimes E)\tto\cong 
H_{\g,K}^p(H_q^0(\Ga,\Sigma,C^\infty(G))\otimes E).
$$
\end{conjecture}

We can prove the higher order version of Borel's result.

\begin{theorem}
For each $q\ge 1$, the inclusions $C_\umg^\infty(G)\hookrightarrow C_\mg^\infty(G)\hookrightarrow C^\infty(G)$ induce isomorphisms in cohomology:
\begin{multline*}
H_{\g,K}^p(H_q^0(\Ga,\Sigma,C_\umg^\infty(G))\otimes E)\tto\cong 
H_{\g,K}^p(H_q^0(\Ga,\Sigma,C_\mg^\infty(G))\otimes E)\\
\tto\cong 
H_{\g,K}^p(H_q^0(\Ga,\Sigma,C^\infty(G))\otimes E).
\end{multline*}
\end{theorem}

\prf
Let $\Omega$ be one of the spaces $\C_\umg^\infty(G)$ or $\C_\mg^\infty(G)$.

We will now leave $\Sigma$ out of the notation.
By Proposition \ref{2.4} we get an exact sequence
$$
0\to H_q^0(\Ga,\Omega)\to  H_{q+1}^0(\Ga,\Omega)
\to H^0(\Ga,\Omega)^{N_{\Ga,\Sigma}(q)}\to 0,
$$
and the corresponding long exact sequences in $(\g,K)$-cohomology.
For each $p\ge 0$ we get a commutative diagram with exact rows
$${\scriptsize
\divide\dgARROWLENGTH by3
\begin{diagram}
\node{H_{\g,K}^p(H_q^0(\Ga,\Omega)\otimes E)}\arrow{e}\arrow{s,l}{\al}
	\node{H_{\g,K}^p(H_{q+1}^0(\Ga,\Omega)\otimes E)}\arrow{e}\arrow{s,l}{\beta}
	\node{H_{\g,K}^p(H^0(\Ga,\Omega)\otimes E)^{N_{\Ga,\Sigma}(q)}}\arrow{s,l}{\ga}\\
\node{H_{\g,K}^p(H_q^0(\Ga,C^\infty(G))\otimes E)}\arrow{e}
	\node{H_{\g,K}^p(H_{q+1}^0(\Ga,C^\infty(G))\otimes E)}\arrow{e}
	\node{H_{\g,K}^p(H^0(\Ga,C^\infty(G))\otimes E)^{N_{\Ga,\Sigma}(q)}.}
\end{diagram}}
$$
Borel has shown that $\ga$ is an isomorphism and that $\al$ is an isomorphism for $q=1$.
We prove that $\beta$ is an isomorphism by induction on $q$.
For the induction step we can assume that $\al$ is an isomorphism.
Since the diagram continues to the left and right with copies of itself where $p$ is replaced by $p-1$ or $p+1$, we can deduce that $\beta$ is an isomorphism by the 5-Lemma.
\qed

By Proposition \ref{2.2} (b) this proof cannot be applied to $\CA(G)$.

{\small Mathematisches Institut\\
Auf der Morgenstelle 10\\
72076 T\"ubingen\\
Germany\\
\tt deitmar@uni-tuebingen.de}

\today

\end{document}